\documentclass[12pt]{amsart}

\usepackage[dvips]{graphicx}
\usepackage[dvips]{color}

\title[Rapidly growing entire functions]{Rapidly growing entire functions
with three singular values}
\author{Sergei Merenkov}
\address{Department of Mathematics\\
University of Illinois\\ 1409 W. Green Street\\ Urbana, IL
61801\\USA}
\thanks{The author was supported by NSF grants DMS-0400636, DMS-0703617, DMS-0244421,
and DMS-0244547.}

\date{\today}

\usepackage{amsmath, amssymb, amsthm, latexsym}
\newcommand\C{{\mathbb C}}
\newcommand\N{{\mathbb N}}

\newcommand\Z{{\mathbb Z}}
\newcommand\R{{\mathbb R}}
\newcommand\dee{\partial}
\newcommand\co{\colon}

\newcommand\OC{{\overline{\mathbb C}}}

\newtheorem{theorem}{Theorem}

\newtheorem{lemma}{Lemma}

\begin{document}


\email{merenkov@math.uiuc.edu}

\abstract{We settle the problem of finding an entire function with
three singular values whose Nevanlinna characteristic dominates an
arbitrarily prescribed function.}
\endabstract

\maketitle

\section{Introduction}\label{S:Intro}

Let $f$ be a transcendental meromorphic function in the plane $\C$. A
\emph{critical point} of $f$ is a point at which the spherical
derivative of $f$ vanishes. The value of $f$ at a critical point
is called a \emph{critical value}. A point $a$ in the sphere $\OC$
is called an \emph{asymptotic value} of $f$ if there exists a
curve $\gamma\co [0,1)\to \C$ such that
$$
\gamma(t)\to\infty\quad {\rm and }\quad f(\gamma(t))\to a\quad
{\rm as }\quad t\to1.
$$
A point $a$ in $\OC$ is a \emph{singular value} of $f$ if it is
either a critical or an asymptotic value. In this paper we study
the growth behavior of entire and meromorphic functions which have
finitely many singular values. The class of such functions is
usually denoted by $\mathcal S$, after A.~Speiser~\cite{aS30},
\cite{aS29}.

%An entire
%function whose Nevanlinna characteristic dominates an arbitrarily
%prescribed function is harder to construct. The reason is that it
%is difficult to arrange a rapid volume growth and still have
%parabolic type in this case.

If $f$ is an arbitrary meromorphic function in the plane, the
Nevanlinna characteristic of $f$ is defined as (see~\cite{wH64},
\cite{Nevbook})
$$
T(r, f)=N(r,f)+m(r,f),
$$
where
$$
N(r,f)=\int_0^r\frac{n(t,f)}{t}\, dt,\ \ \ \ m(r,f)=
\frac1{2\pi}\int_0^{2\pi}\log^+ |f(re^{i\varphi})|\,d\varphi,
$$
and $n(t,f)$ is the number of poles of $f$ in $\{|z|<t\}$. Here we
assumed that 0 is not a pole of $f$. If $f$ is a rational function
of degree $d$, then its Nevanlinna characteristic $T(r, f)$ grows
like $d\log r$, as $r\to\infty$. If $f$ is a transcendental
meromorphic function, then $T(r, f)$ grows faster than any
multiple of $\log r$, but it is easy to see that for any $a>1$ one
can find a transcendental $f$ for which $T(r,f)$ grows slower than
$\log^a r$, as $r\to\infty$.

The question of slowest possible growth of the Nevanlinna
characteristic for meromorphic functions with finitely many
singular values has been studied in recent years, notably by
A.~Eremenko~\cite{aE03}, and J.~Langley~\cite{jL99}, \cite{jL02}.
In particular it was proved that if $f$ is a transcendental
meromorphic function with three singular values, then
$$
\liminf_{r\to\infty}\frac{T(r,f)}{\log^2 r}\geq\frac{\sqrt3}{\pi},
$$
and the constant on the right hand side is sharp. Langley
established the existence of an absolute constant for the right
hand side, and Eremenko found the exact value for this
constant. If $f$ is a transcendental entire function with three
singular values, then $\liminf_{r\to\infty}{T(r,f)}/{\log^2 r}$ is
infinite. In fact, the Nevanlinna characteristic $T(r,f)$ of such
a function dominates a positive multiple of $\sqrt r$.

In general, if $f$ is a transcendental meromorphic function which
has finitely many singular values, then Langley showed that
$$
\liminf_{r\to\infty}\frac{T(r,f)}{\log^2 r}>0,
$$
but the left hand side can be as small as one wishes if the number
of singular values is greater than three.

Here we investigate the question of arbitrarily rapid growth.
%
%\vspace{.2 cm}\noindent \emph{Can one find an entire or a
%meromorphic function with three singular values whose Nevanlinna
%characteristic dominates an arbitrarily prescribed function?}
%\vspace{.2 cm}
%
%The main result of this paper is the following theorem, which
%answers this question in the positive for the stronger case of
%entire functions.
\begin{theorem}\label{T:MT}
For every $\R$-valued function $M(r),\ r\geq0$, there exists an
entire function $f$ with three singular values 0, 1, and $\infty$
such that
$$
T(r, f)\geq M(r),\ {\text{for }}r\geq r_0,
$$
and some $r_0>0$.
\end{theorem}

Our proof of this theorem is based on a combinatorial construction
of a Riemann surface spread over the sphere which branches over
three points. The desired map is obtained as a composition of a
uniformizing map of this Riemann surface and the projection map to
the sphere. One of the key steps in proving Theorem~\ref{T:MT} is
to establish a quantitative control on the volume growth of a
graph in terms of the combinatorial modulus. This is done in
Lemma~\ref{L:KeyL}.

A meromorphic function whose Nevanlinna characteristic dominates
an arbitrarily prescribed function is easier to produce. Indeed,
there is more flexibility in constructing surfaces spread over the
sphere that correspond to meromorphic functions, since one does
not need to worry about $\infty$ being an omitted value. The
construction is outlined in Section~\ref{S:Mero}.

{\sc Acknowledgments.} {The author would like to thank
Alex~Eremenko for suggesting this problem, and Mario~Bonk,
Lukas~Geyer, and Juha~Heinonen for many useful discussions and
interest in this work. Also, many thanks go to the anonymous
referee for numerous useful comments.}

\section{Graphs}

A graph $G$ is a pair $(V_G, E_G)$, where $V_G$ is a set of
\emph{vertices} and $E_G$ is a subset of unordered pairs of
elements in $V_G$, called \emph{edges}. If $v_1, v_2\in V_G$, and
$\{v_1, v_2\}\in E_G$, we say that $v_1$ and $v_2$ are
\emph{connected by an edge}, and write $v_1\sim v_2$. We assume
that no vertex is connected to itself by an edge.
%, and two distinct vertices are not connected by more than one edge.
A graph is called \emph{bipartite} if the vertices can be
subdivided into two disjoint sets, say $A$ and $B$, and every edge
connects a vertex from $A$ to one from $B$. A \emph{subgraph} $G'$
of a graph $G$ is a graph whose vertex set forms a subset of
$V_G$, and if two vertices of $G'$ are connected by an edge in
$G'$, then they are connected by an edge in $G$. If $A$ is a
subset in $V_G$, we denote by $|A|$ the cardinality of $A$, where
$|A|=\infty$ if the set $A$ is infinite.

The {\emph{valence}} of $v\in V_G$ is $|\{u\in V_G\co\ u\sim
v\}|$. The {\emph{valence}} of $G$ is the supremum of the valences
over all vertices of $G$. A graph $G$ is called \emph{locally
finite}, if the valence of each vertex is finite. A graph is said
to have a \emph{finite valence}, if there is a uniform bound on
the valence at each vertex. A graph is called \emph{homogeneous}
of valence $q$ if every vertex has the same valence $q$.

A {\emph{chain}} in $G$ is a sequence $(\dots, x_{-1}, x_0,
x_1\dots)$ of vertices, finite or infinite in one or both
directions, such that $\dots\sim x_{-1}\sim x_0\sim x_1\sim\dots$.
We also refer to a chain as a sequence of vertices along with the
edges connecting them. A chain $(\dots, y_{-1}, y_0, y_1,\dots)$
is a {\emph{subchain}} of a chain $(\dots,x_{-1}, x_0,
x_1,\dots)$, if $y_j=x_{k(j)}$ for some monotone increasing
sequence $(k(j))$. We say that a chain $(x_1,\dots, x_n)$
{\emph{connects}} two subsets $A$ and $B$ of $V_G$, if $x_1\in A$
and $x_n\in B$. A chain $(x_1,x_2,\dots)$ {\emph{connects}} a
finite set $A$ to $\infty$, if $x_1\in A$ and it eventually leaves
every finite set, i.e., for every finite subset $K$ of $V_G$ there
exists $k\in \N$ such that $x_j\notin K$ for $j\geq k$. A set $B$
in $V_G$ is said to \emph{separate} a set $A\subset V_G$ from
$\infty$ if every chain connecting $A$ to $\infty$ has a vertex in
$B$.

A \emph{loop} in a graph is a finite chain $(x_1,\dots, x_n)$ such
that $x_1=x_n$ and all other vertices of the chain are distinct. A
\emph{tree} is a graph that does not contain any loops $(x_1,\dots, x_n)$ with $n>3$. A
\emph{subtree} is a subgraph of a tree.

If $G$ is a graph and $V$ is a subset of the vertex set $V_G$, we
consider the subgraph $G'$ of $G$ \emph{determined by the vertex
set} $V$ to be the graph whose vertex set is $V$, and two vertices
$v_1$ and $v_2$ are connected by an edge in $G'$ if and only if
they are connected by an edge in $G$.

A {\emph{domain}} $D$ in a graph $G$ is a subset of the vertex set
$V_G$ which is {\emph{connected}} in the sense that every two
vertices in $D$ can be connected by a chain all of whose vertices
are in $D$. The \emph{boundary} of $D$ in $G$, denoted by $\dee_G
D$, or $\dee D$ if the graph is understood, is the set of all vertices in $V_G$ that are not in $D$, and
each of which is connected by an edge in $E_G$ to a vertex in $D$.
An \emph{annulus} in a graph $G$ is a subset of $V_G$ whose
complement in $V_G$ consists of two disjoint domains. Not every
graph contains an annulus. A sequence of annuli $(A_k)$ is called
\emph{nested} if the annuli are pairwise disjoint and $A_{k+1}$
separates $A_k$ from $\infty$.

In this paper we only consider \emph{planar} graphs, i.e., graphs
embedded in the plane $\R^2$. If we fix an embedding of a graph
into $\R^2$, then we can speak of \emph{faces} of the graph. These
are complementary components of the image of the graph in the
plane. A \emph{side} of a face is a part of its boundary that is
the image of an edge under the embedding. If $G$ is a planar
graph, one can also define its \emph{dual} $G^*$. The vertices of
$G^*$ are in one to one correspondence with the faces of $G$. Two
vertices of $G^*$ are connected by an edge if and only if the
boundaries of the corresponding faces of $G$ share a side.

A connected graph can be viewed as a metric space if one declares
that every edge is isometric to a unit interval on the real line.
This metric restricts to the space whose elements are vertices of
the graph, in which case it is said that the graph is endowed with
the \emph{word metric}. Thus we can speak of \emph{geodesics} in a
graph, i.e., chains connecting two vertices or two sets and having
the smallest lengths among all such chains. If $A$ and $B$ are two
subsets of $V_G$, we denote by $\delta(A,B)$ the word distance
between $A$ and $B$, i.e., the number of edges in a geodesic
connecting $A$ and $B$. If $A$ is a one vertex set $\{v\}$, we
write $\delta(v,B)$ instead of $\delta(\{v\},B)$. Similarly,
$\delta(v,w)$ stands for $\delta(\{v\},\{w\})$.

\section{Surfaces of Speiser class}

A \emph{surface spread over the sphere} is an equivalence class of
pairs $[(X,\pi)]$, where $X$ is an open, i.e., non-compact, simply
connected topological surface and $\pi\co X\to\OC$ is a
continuous, open, and discrete map. Two pairs $(X_1,\pi_1)$ and
$(X_2,\pi_2)$ are \emph{equivalent} if there exists a
homeomorphism $h\co X_1\to X_2$ such that $\pi_1=\pi_2\circ h$.
The map $\pi$ is called the \emph{projection}.

In a neighborhood of each point $x$ in $X$ the map
$\pi$ is given in some local coordinates (for neighborhoods of $x$ and $\pi(x)$) by $z\mapsto z^k$, where
$k=k(x)\in \N$ is called the \emph{local degree} of $f$ at $x$. If
$k\geq2$, then $x$ is called a \emph{critical point} of $f$, and
in this case the value $f(x)$ is called a \emph{critical value}.
As in the case $X=\C$ and $\pi$ a meromorphic function, $a\in \OC$
is called an \emph{asymptotic value} if there exists a curve
$\gamma\co [0,1)\to X$ such that
$$
\gamma(t)\to\infty\quad {\rm and }\quad \pi(\gamma(t))\to a\quad
{\rm as}\quad t\to1.
$$
Here $\gamma(t)\to\infty$ means that $\gamma(t)$ leaves every
compact set of $X$ as $t\to1$. A point $a$ in $\OC$ is a
\emph{singular value} of $\pi$ if it is either a critical or an
asymptotic value.

According to S.~Sto\"ilow~\cite{sS56}, $X$ supports a complex
structure, the \emph{pullback structure}, in which the map $\pi$
is holomorphic. A surface spread over the sphere is said to have
\emph{parabolic type}, or is called \emph{parabolic}, if $X$
endowed with the pullback structure is conformally equivalent to
the plane. Otherwise it is said to have \emph{hyperbolic type}.
The homeomorphism $h$ in the definition of equivalence is a
conformal map in these pullback structures, and therefore the
conformal type of a surface spread over the sphere is well
defined. For simplicity below we refer to a pair $(X,\pi)$, rather
than an equivalence class, as a surface spread over the sphere.
%, and also call it a \emph{surface}
%unless we want to emphasize the projection map $\pi$.
If $g$ is a uniformizing map for $X$ defined in the complex plane
or the unit disc, then $f=\pi\circ g$ is a meromorphic function.
If $\pi$ omits the value $\infty$, then $f$ is holomorphic. The
surface spread over the sphere $(X,\pi)$ is classically referred
to as the ``surface of $f^{-1}$".

A surface spread over the sphere belongs to \emph{Speiser class}
$\mathcal S$ if $\pi$ has only finitely many singular values. If
$\{a_1,\dots, a_q\}$ is the set of singular values of $\pi$, then
$\pi$ restricted to $\pi^{-1}(\OC\setminus\{a_1,\dots,a_q\})$ is a
covering map. Surfaces spread over the sphere of class $\mathcal
S$ have combinatorial representations in terms of Speiser graphs.

Assuming that $(X,\pi)\in \mathcal S$ and $\pi$ has $q$ singular
values $a_1,\dots, a_q$, we fix an oriented Jordan curve $L$ in 
$\OC$, visiting the points $a_1,\dots, a_q$ in cyclic order of
increasing indices.
%The curve $L$ is usually called a {\emph{base curve}}.
This curve decomposes the sphere into two simply connected
regions.
%$H_1$, the region to the left of $L$, and $H_2$, the
%region to the right of $L$.
Let $L_i$, $i\in\{1,2,\dots,q\}$, be the arc of $L$ from $a_i$ to
$a_{i+1}$ (with indices taken modulo $q$). Let us fix points $p_1$
and $p_2$ in the two complementary components of $L$, and choose
$q$ Jordan arcs $\gamma_1,\dots, \gamma_q$ in ${\OC}$, such that
each arc $\gamma_i$ has $p_1$ and $p_2$ as its endpoints, and has
a unique point of intersection with $L$, which is in $L_i$. These
arcs are chosen to be interiorwise disjoint, that is,
$\gamma_i\cap\gamma_j=\{p_1,p_2\}$ when $i\neq j$. Let $\Gamma'$
denote the graph embedded in $\OC$, whose vertices are $p_1$ and
$p_2$, and whose edges are $\gamma_i,\ i=1,\dots, q$, and let
$\Gamma=\pi^{-1}(\Gamma')$. We identify $\Gamma$ with its image in
$\R^2$ under an orientation-preserving homeomorphism of $X$ onto
$\R^2$.
%Clearly it does not depend on the choice of the points $p_1,\
%p_2$, and the curves $\gamma_i,\ i=1, \dots, q$.
The graph $\Gamma$ is infinite, connected, homogeneous of valence
$q$, and bipartite. The vertices that project to $p_1$ are
labelled $\times$ and the ones that project to $p_2$ are labelled
$\circ$. A graph, properly embedded in the plane and having these
properties is called a {\emph{Speiser graph}}.
%The
%vertices of a Speiser graph $\Gamma$ are traditionally marked by
%$\times$ and $\circ$, such that each edge of $\Gamma$ connects a
%vertex marked $\times$ with a vertex marked $\circ$.
%Each {\emph{face}} of $\Gamma$, i.e., a connected component of
%$\R^2\setminus\Gamma$, has either a finite even number of edges
%along its boundary, in which case it is called an {\emph{algebraic
%elementary region}}, or infinitely many edges, in which case it is
%called a {\emph{logarithmic elementary region}}.
Two Speiser graphs $\Gamma_1, \ \Gamma_2$ are said to be
{\emph{equivalent}}, if there is an orientation-preserving
homeomorphism of the plane which takes $\Gamma_1$ to $\Gamma_2$.
% Below we refer to an equivalence class as a Speiser graph.
Each face of the Speiser graph $\Gamma$ is labelled by the
corresponding element of the set $\{a_1,\dots, a_q\}$.

The above construction of a Speiser graph from a surface spread
over the sphere of class $\mathcal S$ is reversible. Suppose we
are given a Speiser graph $\Gamma$ whose faces are labelled by
$a_1, \dots, a_q$. A necessary condition for existence of a
surface spread over the sphere of class $\mathcal S$ with singular
values $a_1,\dots, a_q$ and whose Speiser graph is $\Gamma$ is
that the labels should satisfy a certain compatibility condition.
Namely, when going counterclockwise around a vertex $\times$, the
indices are encountered in their cyclic order, and around $\circ$
in the reversed cyclic order.  We fix a simple closed curve
$L\subset\OC$ passing through $a_1, \dots, a_q$. Let $H_1,H_2$ be
the complementary regions whose common boundary is $L$, and let
$L_1,\dots,L_{q}$ be as above. Let $\Gamma^*$ be the planar dual
of $\Gamma$. The vertices of $\Gamma^*$ are naturally labelled by
$a_1,\dots, a_q$. If $e$ is an edge of $\Gamma^*$ connecting $a_j$
and $a_{j+1}$, let $\pi$ map $e$ homeomorphically onto the
corresponding arc $L_j$ of $L$. This defines a map $\pi$ on the
edges and vertices of $\Gamma^*$. We then extend $\pi$ to map the
faces of $\Gamma^*$ homeomorphically to $H_1$ or $H_2$, depending
on the orientation of the boundaries. This defines a surface
spread over the sphere $(\R^2, \pi)\in \mathcal S$. The
corresponding labelled Speiser graph is the graph $\Gamma$ with
the prescribed labels. Thus, up to a choice of the curve $L$, we
have a one to one correspondence between surfaces spread over the
sphere of class $\mathcal S$ and equivalence classes of labelled
Speiser graphs. See~\cite{aG70} and \cite{Nevbook} for further details.

\section{Type problem}

A long studied problem is the one of recognizing the conformal
type of a surface spread over the sphere of class $\mathcal S$
from its Speiser graph. An infinite locally finite connected graph
is called \emph{parabolic} if the simple random walk on it is
recurrent. Otherwise it is called \emph{hyperbolic}.
P.~Doyle~\cite{pD84} gave a criterion of type for a surface spread
over the sphere of class $\mathcal S$ in terms of a so called
extended Speiser graph.

Let $\Z_+$ denote the set of non-negative integers. A
{\emph{half-plane lattice}} $\Lambda$ is the graph embedded in
$\R^2$ whose vertices form the set $\Z\times\Z_+$, and $(x',
y')\sim (x'', y'')$ if and only if $(x''-x', y''-y')=(\pm 1, 0)$
or $(0, \pm 1)$. The \emph{boundary} of the half-plane lattice
$\Lambda$ is the infinite connected subgraph of $\Lambda$
determined by the vertex set $\Z\times\{0\}$. There is an action
of $\Z$ on $\Lambda$ by horizontal shifts. A {\emph{half-cylinder
lattice}} $\Lambda_n$ is $\Lambda/n\Z$. The \emph{boundary} of
$\Lambda_n$ is the induced boundary from $\Lambda$.

Suppose that $\Gamma$ is a Speiser graph and let $n\in\N$ be
given. If we replace each face of $\Gamma$ with
$2k$ edges on the boundary, $k\geq n$, by the half-cylinder
lattice $\Lambda_{2k}$, and each face with infinitely many edges
on the boundary by the half-plane lattice $\Lambda$, identifying
the boundaries of the faces with the boundaries of the
corresponding lattices along the edges and vertices, we obtain an
{\emph{extended Speiser graph}} $\Gamma_n$. The graph $\Gamma_n$
is an infinite connected graph embedded in the plane, containing
$\Gamma$ as a subgraph. It has a finite valence, and all faces of
$\Gamma_n$ have no more than $\max\{2(n-1), 4\}$ sides.

\vspace{.2 cm}
%\begin{theorem}\label{T:Ut}~\cite{pD84}
\noindent {\bf Theorem A}~\cite{pD84} A surface spread over the
sphere $(X, \pi)\in \mathcal S$ is parabolic if and only if
$\Gamma_{1}$ is parabolic.
%, where $\Gamma$ is the Speiser graph
%of $(X, \pi)$}.
%\end{theorem}
\vspace{.2 cm}

In~\cite{sM03} we proved a slight modification of Theorem~A.

\vspace{.2 cm} \noindent {\bf Theorem B}~\cite{sM03} Let $n\in\N$
be fixed. A surface spread over the sphere $(X, \pi)\in \mathcal
S$ is parabolic if and only if $\Gamma_{n}$ is parabolic.
%$(X, f)\in \mathcal S$ is hyperbolic, if and only if the
%McKean-Sullivan random walk on its Speiser graph $\Gamma$ is
%transient. In plain terms, the McKean-Sullivan random walk on
%$\Gamma$ comes from the simple random walk on $\Gamma_1$, when we
%observe it only as it hits $\Gamma$.
\vspace{.2 cm}

Doyle's arguments are probabilistic and electrical,
whereas~\cite{sM03} employs geometric methods, using results
of M. Kanai \cite{mK85}, \cite{mK86}. Below we derive Theorem~B from Theorem~A using results from~\cite{DS}.

Kanai shows the invariance of type under quasi-isometries for
spaces with bounded geometry. A map $\Phi: (X_1,d_1)\to (X_2,
d_2)$ between two metric spaces
%, not necessarily continuous,
is called a \emph{quasi-isometry}, if the following conditions are
satisfied:
\begin{enumerate}
    \item[1.] for some $\epsilon>0$, the $\epsilon$-neighborhood
    of the image of $\Phi$ in $X_2$ covers $X_2$;
    \item[2.] there are constants $k\geq 1, \ C\geq 0$, such that
    for all $x_1,\ x_2\in X_1$,
$$
k^{-1}d_1(x_1, x_2) - C\leq d_2(\Phi(x_1), \Phi(x_2))\leq
kd_1(x_1, x_2) + C.
$$
\end{enumerate}
%\end{definition}
The metric space $(X_1, d_1)$ is {\it{quasi-isometric}}
to the metric space $(X_2, d_2)$ if there exists a quasi-isometry
from $X_1$ to $X_2$. This is an equivalence relation. The notion
of quasi-isometry, or rough isometry, was introduced by M.~Gromov
\cite{mG81}.

A Riemannian surface has \emph{bounded geometry} if it is
complete, the Gaussian curvature is bounded from below, and the
radius of injectivity is positive. The latter means that there
exists $\delta>0$ such that every open ball whose radius is at most
$\delta$ is homeomorphic to a Euclidean ball. Kanai proves that
if a Riemannian surface has bounded geometry and is
quasi-isometric to a finite valence graph with the word metric,
then the surface and the graph have the same type. Likewise, two
quasi-isometric graphs with finite valence have the same type.

%Here we deduce Theorem~B from Theorem~A using results
%from~\cite{DS}.

\emph{Proof of Theorem~B.} By Theorem~A one needs to show that
$\Gamma_n$ is parabolic if and only if $\Gamma_1$ is. Assume first
that $\Gamma_1$ is parabolic. The graph $\Gamma_n$ is obtained
from $\Gamma_1$ by cutting the edges that connect the vertices of
$\Gamma$, viewed as a subgraph of $\Gamma_1$ using the obvious
embedding, on the boundary of faces of $\Gamma$ with $2k$ edges, $k<n$, to the
vertices of $\Lambda_{2k}$. Therefore this direction follows from
the Cutting Law~\cite{DS}, p.100. For the other direction, assume that
$\Gamma_1$ is hyperbolic. We consider a new graph
$\tilde\Gamma_1$, obtained from $\Gamma_1$ by shorting all
non-boundary vertices of every half-cylinder lattice
$\Lambda_{2k},\ k<n$, that have replaced a face of $\Gamma$. Here
\emph{shorting} a set of vertices means identifying them. By the
Shorting Law~\cite{DS}, p.100, $\tilde\Gamma_1$ is also
hyperbolic. But $\tilde\Gamma_1$ has finite valence
and is quasi-isometric to $\Gamma_n$. The quasi-isometry is
given by an embedding of $\Gamma_n$ into $\tilde\Gamma_1$ induced
from the obvious embedding of $\Gamma_n$ into $\Gamma_1$.
Therefore $\Gamma_n$ is hyperbolic. \qed

Due to the nature of a construction, as in our case below, it is
often easier to establish the type for the dual graph $\Gamma_n^*$
to the extended Speiser graph $\Gamma_n$.

\vspace{.2 cm} \noindent {\bf Theorem C } Let $n\in\N$ be fixed. A
surface spread over the sphere $(X,\pi)\in\mathcal S$ is parabolic
if and only if $\Gamma_n^*$ is parabolic.

\vspace{.2 cm} \emph{Proof.} The graph $\Gamma_n$ in question and
its dual have finite valence. A map $\Phi$ that sends every
vertex $v$ of $\Gamma_n^*$ to any vertex on the boundary of the
face of $\Gamma_n$ corresponding to $v$ is a quasi-isometry.
Indeed, the first condition for quasi-isometry follows since every
vertex of $\Gamma_n$ is on the boundary of a face and there is a
uniform bound on the number of sides of each face since
$\Gamma_n^*$ has finite valence. Therefore every vertex of
$\Gamma_n$ is within a uniformly bounded distance from an image of
a vertex in $\Gamma_n^*$ under $\Phi$.

The second condition follows since both graphs have finite
valence. Let $\gamma^*$ be a geodesic chain in $\Gamma_n^*$
connecting two vertices $v_1$ and $v_2$. By tracing the boundaries
of faces corresponding to the vertices of $\gamma^*$, one can find
a chain in $\Gamma_n$ connecting $\Phi(v_1)$ and $\Phi(v_2)$, and
whose length is at most $C_1$ times the length of $\gamma$, where
$C_1$ depends only on the valences of $\Gamma_n$ and $\Gamma_n^*$.
Conversely, for every geodesic chain $\gamma$ in $\Gamma_n$
connecting two vertices $\Phi(v_1)$ and $\Phi(v_2)$, one can find
a chain in $\Gamma_n^*$ that connects $v_1$ and $v_2$ by following
the faces that contain $\gamma$ on their boundaries, such that the
length of this new chain is at most $C_2$ times the length of
$\gamma$. The constant $C_2$ depends only on the valences of
$\Gamma_n$ and $\Gamma_n^*$.

Since the graphs $\Gamma_n$ and $\Gamma_n^*$ are quasi-isometric
and have finite valences, they have the same type. Now the result
follows from Theorem~B.\qed

\section{Combinatorial modulus}\label{S:Cm}

In 1962 R.~J.~Duffin~\cite{rD62} introduced a combinatorial
modulus for chain families in graphs. In his setting the masses
are assigned to the edges of the graph. Parabolicity of a locally
finite graph is equivalent to the condition that the modulus of
the family of chains connecting a fixed vertex to infinity is
zero. For our purposes it is more convenient to use a different
notion of modulus, introduced more recently by
J.~W.~Cannon~\cite{jC94}, where masses are assigned to vertices
rather than edges. This approach leads to certain combinatorial
uniformization results, see e.g.~\cite{oS93}. If a graph has finite valence, as in our case below, it does not
matter which definition of combinatorial modulus one uses when
establishing parabolicity. This can be seen by distributing masses
from vertices to edges and \emph{vice versa}.

A {\emph{mass distribution}} for a graph $G$ is a non-negative
function on $V_G$. Let $\mathcal C$ be a family of chains in $G$.
We say that a mass distribution $m$ is \emph{admissible} for
$\mathcal C$, if for each chain $(\dots, x_{-1}, x_0,
x_1,\dots)\in\mathcal C$, its \emph{weighted length} $\sum
m(x_j)\geq1$. We denote by ${\rm{mod }}_G\mathcal C$ the
{\emph{combinatorial modulus}} of the chain family $\mathcal C$,
namely
$$
{\rm{mod }}_G\mathcal C=\inf\bigg\{\sum m(v)^2\bigg\},
$$
where the infimum is taken with respect to all admissible mass
distributions, and the sum is over all vertices in $V_G$. We write
${\rm{mod}}\,\mathcal C$ if the graph is understood. To
distinguish, the conformal modulus of a curve family on a surface
will be denoted by $\rm{Mod }$. If $\mathcal C$ is the family of
all chains connecting sets $A$ and $B$, or a set $A$ to $\infty$,
we denote ${\rm{mod}}\,\mathcal C$ by ${\rm{mod }}(A, B)$ or
${\rm{mod }}(A,\infty)$, respectively. If $A$ is an annulus in a
graph $G$, then ${\rm mod }A$ denotes the modulus of the family of
all chains that connect the complementary components of $A$ in
$V_G$.

As for the classical conformal modulus, if $\mathcal C$ and
$\mathcal C'$ are two families of chains, such that every chain in
$\mathcal C$ contains a subchain which is in $\mathcal C'$, then
${\rm{mod}}\,\mathcal C\leq{\rm{mod}}\,\mathcal C'$. Also, if
$(A_k)$ is a sequence of (disjoint) nested annuli, then
$$
{\rm mod }(\{v_0\},\infty)\leq\frac1{\sum1/{\rm mod }A_k},
$$
for any vertex $v_0$ that is separated from $\infty$ by every
$A_k$. The first property follows immediately from the definition.
A proof of the inequality mimics that for the classical conformal
modulus.
%If a graph has finite valence, as in our case below, it does not
%matter which definition of combinatorial modulus one uses when
%establishing parabolicity via the criterion stated in the first
%paragraph of this section. This can be seen by distributing masses
%from vertices to edges and \emph{vice versa}. Thus t
Now, as in the classical case, to show
parabolicity of a finite valence graph, it is enough to exhibit a
sequence $(A_k)$ of (disjoint) nested annuli, such that
$$
\sum\frac1{{\rm mod }A_k}=\infty.
$$
This will be used in the proof of Lemma~\ref{L:KeyL}.

\section{Meromorphic example}\label{S:Mero}

Since later we prove that the Nevanlinna characteristic of an
entire function dominates an arbitrarily prescribed function, here
we only give an outline that such a meromorphic function exists.

Consider the infinite half-strip in the plane
$$
S=\{z=x+iy\co 0\leq x\leq 2,\ 0\leq y< \infty\},
$$
subdivided into squares
$$
\{z\co j\leq x\leq j+1,\ n\leq y\leq n+1\},\quad j= 0,1,\quad n=0,
1, 2,\dots.
$$
For each even $n=2k$, we attach $N(k)$ Euclidean squares with side
length 1 to the edge
$$
e_k=\{z\co x=1,\ 2k\leq y\leq 2k+1\},
$$
so that all of these squares share the side $e_k$, and are
otherwise disjoint. More specifically, we cut the strip $S$ along
$e_k$, take a two-sided unit square cut along one of its edges,
and glue the square to the strip along a cut. We repeat this
operation if necessary, attaching more squares to $e_k$. What
results can be thought of as a book spread open along its spine.
%One can
%think of these attached squares as ``sticking out" from the plane.

The result of the gluing of all the squares is a simply connected
Riemann surface $Y$ with boundary, which corresponds to the
boundary of $S$. Now we consider the double $X$ of $Y$ across the
boundary. This means that $X$ is obtained from two copies of $Y$ by identifying every boundary point of one copy with the point of the other copy that corresponds to the same point of $Y$. This is a simply connected Riemann surface without
boundary. For each $n=0,1,2,\dots$, let $A_n$ denote the annulus
in $X$ that consists of all points corresponding to the points of
the horizontal rectangle $\{n\leq y\leq n+1\}$ of $S$ and all points
of squares attached to $e_{n/2}$ if $n$ is even.
%, and all the corresponding points of the double. 
Each surface $X$ is parabolic
since it contains a sequence of annuli $(A_n)$, where $n$ is odd,
of fixed modulus. Using a modulus estimate, one can show that if
$F$ is a uniformizing map of $\C$ onto $X$, then the image $I_r$
under $F$ of the disc $D_r$ centered at 0 of radius $r$ contains a
ball (in the intrinsic metric of $X$) of radius
$$
L(r)\geq C\log r,
$$
where $C$ is a constant not depending on the sequence $(N(k))$.
Indeed, let $s$ denote the set in $X$ that corresponds to the
segment in $S$ connecting $(0,0)$ to $(2,0)$, and let $s_F$ be the
preimage of $s$ under $F$. The set $s_F$ is homeomorphic to a line segment. Suppose that $n(r)$ is the smallest
natural number so that the annulus $A_{n(r)}$ is not contained in
$I_r$. The conformal modulus of the curve family consisting of
curves in $D_r$ that separate $s_F$ from the boundary of $D_r$
grows like $\log r/(2\pi)$ as $r\to\infty$. On the other hand, the
conformal modulus of the image family in $X$ is bounded above by
$C'n(r)$, where $C'$ is a constant independent of $(N(k))$. This
can be seen by choosing a weight function equal 1/2 at all points
of the annuli $A_0, A_1,\dots, A_{n(r)+1}$ that correspond to
points of $S$, and equal 0 at all other points of these
annuli. From the invariance of modulus under conformal maps we
obtain that
$$
n(r)\geq \log r/(2\pi C'),
$$
which immediately gives the desired estimate for $L(r)$.

Now, by choosing $N(k)$ to grow sufficiently rapidly, one can
arrange arbitrarily rapid growth of the areas, with respect to the
radii, of the intrinsic balls of $X$ centered at some point.
Arbitrarily rapid growth of the areas implies arbitrarily rapid
growth of the Nevanlinna characteristic (see Ahlfors-Shimizu
characteristic in~\cite{wH64}). A similar fact is based on the First Main Theorem of Nevanlinna and it will be discussed in Section~\ref{S:VG}.

By subdividing each square of the surface $X$ into four triangles
using diagonals, and considering the Speiser graph which is dual
to such a triangulation, we obtain a meromorphic function with
three singular values that has the desired properties.

\section{Entire functions with three singular values}

If $f$ is a transcendental entire function with three singular
values $0, 1$, and $\infty$, then $f^{-1}([0, 1])$ forms a locally
finite, infinite tree $T$ embedded in $\R^2$. The
vertices are the preimages of 0 and 1, and the edges are the
preimages of $[0,1]$. Indeed, the graph is connected since $f$
restricted to $f^{-1}(\C\setminus\{0,1\})$ is a covering map. The
valence of each vertex is the local degree of $f$ at the
corresponding point. The graph is infinite since $f$ is
transcendental. Finally, it is a tree because otherwise there
would exist a complementary component of $f^{-1}([0,1])$ that is
compactly contained in $\C$. This is impossible since such a
component would have to contain a preimage of $\infty$, but $f$ is
assumed to be entire.

Conversely, suppose we are given an arbitrary locally finite,
infinite, embedded tree $T$, whose vertices are labelled 0 and 1,
and each edge connects 0 and 1. We construct a surface spread over
the sphere $(X,\pi)$ with three singular values as follows. For
every vertex $v$ in $V_T$ of valence $k$ we consider $k$
non-homotopic, non-intersecting Jordan arcs in $\R^2\setminus T$
that originate at $v$ and escape to infinity. We can choose the
arcs corresponding to different vertices to be disjoint. This
gives a triangulation $T'$ of $\R^2$, with each triangle having an
ideal vertex at infinity. Every triangle of $T'$ has an edge of
$T$ and two arcs escaping to infinity as its sides. Each vertex of
$T'$ has an even valence, and it receives a label 0 or 1 from the
corresponding label of $T$. The ideal vertices at infinity are
labelled by $\infty$.

Consider the dual graph to $T'$, denoted $\Gamma$. The graph
$\Gamma$ is an infinite connected graph, properly embedded in the
plane. It has valence three at each vertex, and every face of
$\Gamma$ has an even (or infinite) number of vertices on its
boundary, so $\Gamma$ is bipartite. Therefore $\Gamma$ is a
Speiser graph. Let $(X,\pi)$ denote a surface spread over the
sphere that corresponds to $\Gamma$ with the induced labels from
$T'$, which are 0, 1, and $\infty$. These are the singular values
of $\pi$, and $\pi$ omits the value $\infty$. Thus the composition
of a uniformizing map of $X$ with $\pi$ is a holomorphic function.
We proceed by explicitly describing $(X,\pi)$ up to conformal
equivalence.
%This has
%the advantage of being more visual and we found it easier to
%describe the example of a desired entire function using a ``dual"
%construction rather than using a Speiser graph.

Let
$$
\alpha=\{(x,y)\co 0\leq x,\ 0\leq y\leq 1\}
$$
be a half-strip in the plane. To each triangle $t$ of $T'$ we
associate a copy of $\alpha$, which we denote by $\alpha(t)$, so
that under an orientation-preserving homeomorphism of the plane
the side of $t$ contained in $T$ corresponds to the segment
joining $(0,0)$ and $(0,1)$, and the sides of $t$ that are in $T'\setminus T$ correspond to two horisontal
rays. If $t_1$ and $t_2$ are adjacent triangles, we glue
$\alpha(t_1)$ and $\alpha(t_2)$ along the corresponding sides
using the identity map. The result of the gluing is a simply
connected open Riemann surface, which we denote by $S(T)$. A tree isomorphic to $T$ embeds in $S(T)$, and we identify this tree with $T$.
%whose edges are images of the vertical sides of $\alpha(t)$ under
%the gluing maps.
%For convenience, we also denote this tree by $T$.
Now we consider the conformal map, continuously extended to the boundary, from the half-strip
$$
\alpha^o=\{(x,y)\co 0<x,\ 0<y<1\}
$$
to the lower half-plane that takes $(0,0), (0,1)$ and $\infty$ to
0, 1, and $\infty$, respectively. This map extends by reflection to a
conformal map from the Riemann surface $S(T)$ to the surface
spread over the sphere $(X,\pi)$ with the pullback complex
structure. The tree $T$ is isomorphic to $\pi^{-1}([0,1])$ with
the natural graph structure.

Since we need to consider an extended Speiser graph in deciding
the type of a surface spread over the sphere, the following
subdivision of $S(T)$ is useful. We subdivide $\alpha$ into squares
$$
\alpha_k=\{(x+k, y)\co 0\leq x\leq1,\ 0\leq y\leq 1 \},\
k=0,1,2,\dots.
$$
The subdivision of $\alpha$ by $\alpha_k,\ k=0,1,2,\dots$, induces
a subdivision of $S(T)$ into squares, a \emph{square subdivision}.
The 1-skeleton of this subdivision considered as a graph will be
denoted by $\sigma=\sigma(T)$. The tree $T$ is a subgraph of
$\sigma$. In the case when the tree $T$ has valence $n$, as in our
example below with $n=3$, the graph $\sigma(T)$ is the dual graph
$\Gamma_n^*$ to the extended Speiser graph $\Gamma_n$. According
to Theorem~C, the surface spread over the sphere $(X,\pi)$ is
parabolic if and only if $\sigma$ is.

\section{Volume growth}\label{S:VG}

The First Main Theorem of Nevanlinna (see~\cite{wH64},
\cite{Nevbook}) asserts that for every $a\in\C$,
$$
T(r,f)=N\bigg(r,\frac1{f-a}\bigg)+m\bigg(r,\frac1{f-a}\bigg)+O(1),\
r\rightarrow\infty.
$$
Therefore, by choosing $a$ to be either 0 or 1, we conclude that
in order to find $f$ with $T(r,f)$ growing arbitrarily rapidly, it
is sufficient to find an embedded tree $T$ with the following
properties. The corresponding surface $S(T)$ is parabolic, and if
$M(r),\ r\geq0$, is a prescribed function, and $g$ a uniformizing
map from $\C$ to $S(T)$, then the number of vertices of
$g^{-1}(T)$ in the disc of radius $r$ about 0 is greater than
$M(r)$, for all $r\geq r_0>0$. In this case the first term
$N(r,1/(f-a))$ alone dominates $M(r)$.

Assuming that $S(T)$ is parabolic and $g$ is a uniformizing map
from $\C$ to $S(T)$,  we denote by $n(r,T,g)$ the number of
vertices of $g^{-1}(T)$ contained in the disc of radius $r$
centered at 0. This is an analog of the counting function $n(r,f)$
in the definition of Nevanlinna characteristic $T(r,f)$.
Theorem~\ref{T:MT} follows from the following theorem, proved in
Section~\ref{S:Pr}.
\begin{theorem}\label{T:main}
Given any $\R$-valued function $M(r),\ r\geq0$, there exists a
locally finite, infinite tree $T$, embedded in the plane, such
that $S(T)$ is parabolic, and $n(r,T,g)\geq M(r), \ r\geq r_0$,
for any uniformizing map $g$ and  $r_0=r_0(g)>0$. Moreover, we can
choose $T$ to be a subtree of the regular tree of valence
three, denoted $T_3$.
\end{theorem}

The tree $T_3$ is homogeneous of valence 3, and we think of $T_3$
as being embedded in the plane. Let $v_0$ be a fixed vertex in
$V_{T_3}$, and $\epsilon_0$ denote the combinatorial modulus
${\rm{mod }}_{T_3}(\{v_0\},\infty)$, which is a positive number
because $T_3$ is hyperbolic, as is well-known. The complement of
$T_3$ in the plane has infinitely many components, three of which
have $v_0$ on their boundaries. We consider one of
these three components, denoted $D$, and let $c=(\dots, v_{-1},
v_0, v_1,\dots)$ be the chain in $T_3$ such that $v_j\neq v_k$ for
$j\neq k$, and $c$ together with the edges that connect its
vertices bounds $D$.

If $k\in\N$, then $T_3\setminus\{v_k,v_{-k}\}$ is a
union of five disjoint domains, one of which contains $v_0$, and each of the four others is bounded by either $v_k$ or $v_{-k}$.  
For each $k\in\N\cup\{0\}$, let $\mathcal C_k$ be the family of
all chains $(x_1,x_2,\dots)$ in $T_3$ that connect $\{v_0\}$ to
$\infty$, and such that all but finitely many of $x_j$'s are
contained in one of the domains into which
$T_3\setminus\{v_k,v_{-k}\}$ splits, that does not contain $v_0$.
The family $\mathcal C_0$ consists of all chains connecting $v_0$
to $\infty$. If $k>0$, each chain of $\mathcal C_k$ should have all but finitely many of its vertices to lie in one of the four domains of $T_3\setminus\{v_k,v_{-k}\}$ that does not contain $v_0$. In other
words, every chain in $\mathcal C_k$ should escape to infinity
through either $v_k$ or $v_{-k}$. It is easy to see that the
sequence $(\epsilon_k)$ defined by $\epsilon_k={\rm{mod }}_{T_3}
\mathcal C_k$ decreases, $0<\epsilon_k\leq\epsilon_0$ for every
$k\in\N$, and $\lim\epsilon_k=0$.

For two quantities $a$ and $b$ we use the notation $a\lesssim b$ if
there exists a constant $C>0$ which depends only on the data of
an underlying space, such that $a\leq C b$.
The key step in the proof of
Theorem~\ref{T:main} is the following lemma.
\begin{lemma}\label{L:KeyL}
Let $c, \mathcal C_k$, and $\epsilon_k$ be as above, $k\in
\N\cup\{0\}$. Let $L(\epsilon),\ 0<\epsilon\leq\epsilon_0$, be a
positive decreasing function, $L(\epsilon_0)\geq1$. Let $B_k'$ be
the subset of vertices of $T_3$ defined by
$$
B_k'=\{v\in V_{T_3}\co\delta(v,v_0)=\delta(v,c)+k\},
$$
and let $B_k$ be the subset of $B_k'$ given by
$$
B_k=\{v\in B_k'\co\delta(v,c)\leq L(\epsilon_{k+1})\}.
$$
Then the subtree $T$ of $T_3$, determined by the vertex set
$$
V_T=\cup_{k=0}^{\infty}B_k,
$$
satisfies the property that for every
$\epsilon\in(0, \epsilon_0]$, and every domain $D$ in $T$
with $v_0\in D$, we have
\begin{equation}\label{E:Me}
{\rm{mod }}_T(\{v_0\},\dee D)<\epsilon\ \Rightarrow\
|D|>L(\epsilon).
\end{equation}
Moreover, if
\begin{equation}\label{E:Pa}
2^{[L(\epsilon_1)]}+\dots+2^{[L(\epsilon_k)]}\leq C
2^{[L(\epsilon_{k+1})]},\ \ k=1,2,\dots,
\end{equation}
where $C$ is a positive constant, then $S(T)$ is parabolic.
\end{lemma}
{\emph{Proof.}} It follows from the definition that $B_k',\
k=0,1,2,\dots$, are disjoint, $\cup_{k=0}^{\infty}B_k'=V_{T_3}$,
and every chain in $\mathcal C_k$ has all but finitely many of its
vertices in $\cup_{l\geq k}B_l'$.

Suppose that $\epsilon\in(0,\epsilon_0]$, and let $D$ be a domain
in $T$ with $v_0\in D$, and such that ${\rm{mod }}_T(\{v_0\}, \dee
D)<\epsilon$. There exists $k\in\N\cup\{0\}$ such that
$\epsilon_{k+1}<\epsilon\leq\epsilon_k$. Assume for contradiction
that $|D|\leq L(\epsilon)$. Since $L$ is decreasing, $|D|\leq
L(\epsilon_{k+1})$, and therefore every chain in $\mathcal C_k$
contains a subchain in $T$ that connects $\{v_0\}$ to $\dee_TD$,
the boundary of $D$ in $T$. Indeed, $D$ can also be
considered as a domain in $T_3$, and it cannot contain vertices of
$B_l',\ l\geq k$, that are more than distance
$[L(\epsilon_{k+1})]-1$ away from $c$ because $|D|\leq
L(\epsilon_{k+1})$. Thus every boundary vertex of $D\subset T_3$
contained in $\cup_{l\geq k}B_l'$ is a boundary vertex of
$D\subset T$. Since $v_0\in D$, every chain $c'$ in $\mathcal C_k$ has a subchain
connecting $v_0$ to some boundary vertex $v'$ of $D$ in $T_3$.
Furthermore, $c'$ contains a subchain connecting $v_0$ to $v'\in\dee_TD$. If not, let $v''$ be
the last vertex of $c'$ that belongs to the boundary of $D$ in
$T_3$. Since $D$ is a domain, and hence is connected, and $T_3$ is
a tree, $v''$ either belongs to $c$ or is contained in $\cup_{l\geq
k}B_l'$. But $c$ is contained in $T$, and in the latter case $v''$
belongs to $T$ as a boundary vertex of $D\subset T_3$ contained in
$\cup_{l\geq k}B_l'$. The desired subchain is obtained by removing
edges of $c'$ that connect vertices outside of $V_T$.

Now we have
$$
\epsilon_k={\rm{mod}}_{T_3}\,\mathcal C_k\leq{\rm{mod
}}_{T}(\{v_0\},\dee D)<\epsilon.
$$
This last estimate contradicts our understanding that $\epsilon\leq\epsilon_k$, and proves~(\ref{E:Me}).

%\begin{figure}[htbp]
%\begin{center}
%{\input{sigma.pstex_t}} \caption{ Graph $\sigma(T)$, tree $T$
%(bold lines), and annulus $A_2$ (bold vertices). } \label{f.sigma}
%\end{center}
%\end{figure}

%%\begin{figure}
%%\centerline{\AffixLabels{
%%\includegraphics*[height=2.7in]{sigma.eps}
%%}} \caption{\label{f.sigma}Graph $\sigma(T)$, tree $T$ (bold), and
%%annulus $A_2$ (bold vertices).}
%%\end{figure}

It remains to prove that under Condition~(\ref{E:Pa}), $S(T)$ is
parabolic. The tree $T$ has an axis of symmetry passing through
$v_0$ so that  under the symmetry transformation the
vertex $v_k$ is mapped to $v_{-k}$ and \emph{vice versa}, and each
$B_k$ as well as the chain $c$ are invariant. One should think of this axis of symmetry as being orthogonal to $c$. Let
$\sigma=\sigma(T)$ be the 1-skeleton of the square subdivision of
$S(T)$ that was created using the $\alpha_k$'s. The graph $\sigma$
has also an axis of symmetry, denoted $a$, induced by the
axis of symmetry of $T$. We claim that $\sigma$ is
parabolic. For that purpose we exhibit a sequence of nested annuli
$(A_k)$ and verify that $\sum1/{\rm mod }A_k=\infty$.
%(see
%Figure~\ref{f.sigma}).
%Since the valence of $\sigma$ is finite,
%A result of P.~Doyle~\cite{pD84} (see also~\cite{sM03} for an
%alternative proof) implies that $S(T)$ is parabolic if and only if
%$\sigma$ is parabolic, meaning that ${\rm{mod
%}}_{\sigma}(\{v_0\},\infty)=0$.
%To show this it is enough to
%exhibit a sequence $(A_k)$ of nested annuli separating the
%vertex $v_0$ from $\infty$ such that if ${\rm mod }A_k$ denote the
%connecting

For each $k=1,2,\dots$, we consider an annulus $A_k$ in $\sigma$
obtained as follows.
%(see Figure~\ref{f.sigma}).
The vertices of $T$ separate those of
$\sigma$ into two groups, which we call $V_+$ and $V_-$. The
sets $V_+$ and $V_-$ form the sets of vertices of the upper and
lower square grids $\{(m,n)\co m\in \Z,\ n\in\N\}$ and
$\{(m,n)\co\ m\in\Z,\ -n\in\N\}$, respectively, so that for each
of these sets the vertices with coordinates $(0,n)$ are located on
the symmetry axis $a$.
%The point $(0,0)$
%on the boundary of $V_+$ corresponds to $v_0$, and $(0,0)$ on the
%boundary of $V_-$ corresponds to either $v_0$ if
%$L(\epsilon_1)<1$, or the vertex $v_0'$ in $V_T$ such that
%$\delta(v_0', v_0)=1$, and $v_0'\neq v_1, v_{-1}$.
Each $A_k$ consists of the vertices of the set $B_k\subset V_T$
defined above, vertices $(m,n)$ in $V_+$ such that
$\max\{|m|,|n|\}=k$, and vertices $(m,n)$ in $V_-$ such that
$a_k\leq\max\{|m|, |n|\}\leq b_k$, where $a_k$ and $b_k$ are
chosen as follows. The number $a_k$ is the least one such that
the vertex $(a_k,-1)$ of $V_-$ is connected by an edge to $v_k$,
and $b_k$ is the largest number such that $(b_k,-1)$ is connected
by an edge to $v_k$. A direct calculation gives
$$
%\begin{aligned}
a_k=2^{[L(\epsilon_1)]}+2\bigg(2^{[L(\epsilon_2)]}+\dots+
2^{[L(\epsilon_k)]}\bigg)-k+1,
$$
$$
b_k=2^{[L(\epsilon_1)]}+2\bigg(2^{[L(\epsilon_2)]}+\dots+
2^{[L(\epsilon_{k+1})]}\bigg)-k-1.
%\end{aligned}
$$
Indeed, for each $l>0$, the number of vertices of $B_l$ lying to
one side of the axis of symmetry $a$ is $2^{[L(\epsilon_{l+1})]}$,
and the total number of vertices $v$ of $V_-$ to one side of $a$,
such that $v$ is connected to a vertex in $B_l$, is
$2^{[L(\epsilon_{l+1})]+1}-1$. Adding the latter terms for
$l=1,2,\dots, k-1$ and for $l=1,2,\dots, k$ together, each along
with $2^{[L(\epsilon_1)]}$, contributed by $B_0$, we obtain the
quantities $a_k$ and $b_k+1$, respectively.

Now we assign mass 1 to all vertices in $A_k\cap V_+$, mass
$1/2^{l-1}$ to vertices $v$ in $B_k$ such that $\delta(v,c)=l,\
l=1,2,\dots, [L(\epsilon_{k+1})]$, and mass
$1/2^{[L(\epsilon_{k+1})]-1}$ to the vertices in $A_k\cap V_-$.
This is an admissible mass distribution for the family of chains
that connect the two components of $V_{\sigma}\setminus A_k$.
Indeed, if a chain contains a vertex in $A_k\cap V_+$, we are
done. If a chain only contains vertices of $A_k\cap V_-$, then its
weighted length is at least
$$
\frac{b_k-a_k}{2^{[L(\epsilon_{k+1})]-1}}=
4\bigg(1-\frac1{2^{[L(\epsilon_{k+1})]}}\bigg)\geq1,
$$
since we assumed that $L(\epsilon_0)\geq1$, and $L$ is decreasing.
A chain that contains only vertices of $B_k$ has weighted length
at least 1, because the subgraph of $\sigma$ determined by the
vertex set $B_k$ is a tree, and hence such a chain has to contain
the vertex $v_k$. The remaining case is when a chain $\gamma$
contains vertices of $A_k\cap V_-$ as well as vertices in $B_k$.
It is easy to see that then there is a chain that contains only
vertices of $A_k\cap V_-$, and whose weighted length is comparable
to that of $\gamma$, with absolute constants. Such a chain is
obtained by replacing each vertex $v$ of $\gamma$ that belongs to
$B_k$ by a chain of vertices in $A_k\cap V_-$ of the form
$(m,-1)$, so that the first and the last vertices of this chain
are connected by edges in $\sigma$ to $v$. Multiplying the mass
distribution by an appropriate constant produces an admissible
mass distribution.

The mass bound is
$$
\begin{aligned}
&\lesssim
k+\sum_{l=1}^{[L(\epsilon_{k+1})]}\frac{2^l}{2^{2(l-1)}}+\frac{(2b_k)^2-(2a_k)^2}
{2^{2([L(\epsilon_{k+1})]-1)}}\\
%\frac{(2^{[L(\epsilon_1)]}+\dots+2^{[L(\epsilon_{k+1})]})^2-
%(2^{[L(\epsilon_1)]}+\dots+2^{[L(\epsilon_k)]})^2}{2^{2([L(\epsilon_{k+1})]-1)}}\\
&\lesssim k+1+\bigg(1+2\frac{2^{[L(\epsilon_1)]}
+\dots+2^{[L(\epsilon_k)]}}{2^{[L(\epsilon_{k+1})]}}\bigg)
\lesssim k,\ \ k=1,2,\dots.
\end{aligned}
$$
Since $\sum^{\infty}1/k=\infty$, we conclude that $\sigma$
is parabolic.  \qed

\section{Comparison of moduli}\label{S:Mc}

The results of this section are essentially contained
in~\cite{BK02}, Section~8.

A pathwise connected metric measure space $(X, d, \mu)$ is an $n$-\emph{Loewner space} if
$$
\inf\{{\rm Mod}_n(E,F)\co \Delta(E, F)\leq t\}
$$
is a positive function for all $t>0$, where ${\rm Mod}_n(E,F)$
denotes the $n$-modulus of a curve family connecting two disjoint
continua $E$ and $F$ in $X$, and
$$
\Delta(E,F)=\frac{{\rm dist}(E,F)}{\min\{{\rm diam} E, {\rm diam}
F\}}
$$
is called the \emph{relative distance} between $E$ and $F$.
Loewner spaces were introduced in~\cite{HK98}, see
also~\cite{jH01}.

Recall that $\sigma=\sigma(T)$ is the 1-skeleton of the square
subdivision of $S(T)$. Let $\mathcal U=\{U_v\co\ v\in
V_{\sigma}\}$ be an open cover of $S(T)$, where $U_v$ is the
interior of the union of all squares in $\sigma$ that have a
vertex at $v\in V_{\sigma}$. If $J>0$, we define the
\emph{$J$-star} of $v\in V_{\sigma}$ as
$$
St_J(v)=\cup\{U_u\co u\in V_{\sigma},\ \delta(u, v)< J\}.
$$
Note that $St_1(v)=U_v$. Since $T$ is a tree, it is easy to see
that $St_J(v)$ is an open, connected, and simply connected subset
of $S(T)$. For a set $A$ in $S(T)$ we denote by $V_A$ the set of
vertices $v$ such that $U_v\cap A\neq\emptyset$.

\begin{lemma}\label{L:Ad}
Assume that the valence $k$ of $T$ is finite. Let $v$ be a vertex
of $\sigma$, and $\rho$ be an arbitrary Borel measurable
non-negative function on $St_2(v)$. If $Y_1,\ Y_2\subset S(T)$ are
continua with $Y_i\cap U_v\neq\emptyset$, and ${\rm{diam
}}(Y_i)\geq c_0>0,\ i=1,2$, then there is a rectifiable curve
$\eta$ in $St_2(v)$ connecting $Y_1$ and $Y_2$ such that
$$
\int_{\eta}\rho\, ds\leq
C_0\bigg(\int_{St_{2}(v)}\rho^2d\mu\bigg)^{1/2},
$$
where $C_0>0$ depends only on $c_0$ and $k$.
\end{lemma}
\emph{Proof.} The result follows from the observation that there
are only finitely many, depending on $k$, different possibilities
for $St_2(v)$ that can occur, and from Theorem~{6.13}
in~\cite{HK98}, which implies that $St_2(v)$ is a 2-Loewner space.
Indeed, the Loewner property gives that the conformal modulus
${\rm{Mod }}(Y_1, Y_2)\geq c>0$, where $c$ depends on $c_0$ and
$k$ only. This means that for every Borel measurable non-negative
function $\rho$ on $St_2(v)$ we have
$$
\int_{St_2(v)}\rho^2d\mu\geq
c\inf_{\gamma}\bigg(\int_{\gamma}\rho\, ds\bigg)^2,
$$
where the infimum is taken over all curves $\gamma$ in $St_2(v)$
that connect $Y_1$ and $Y_2$. Thus, for every $\epsilon>0$ there
exists a rectifiable curve $\eta\subset St_2(v)$ connecting $Y_1$
and $Y_2$ such that
$$
\bigg(\int_{\eta}\rho\,
ds\bigg)^2\leq\frac1{c}\int_{St_2(v)}\rho^2d\mu+\epsilon.
$$
Choosing $\epsilon=\frac1{c}\int_{St_2(v)}\rho^2d\mu$ completes
the proof in the case when $\rho$ is not zero almost everywhere on
$St_2(v)$. The latter case is trivial. \qed

\begin{lemma}\label{L:Mm}
If $T$ is an infinite embedded tree of valence $k$, then there
exists a constant $C_1\geq1$, depending only on $k$, such that if
$A, B\subset S(T)$ are two continua not contained in any set
$St_{2}(v)$ for $v$ a vertex of $\sigma$, then
\begin{equation}\label{E:Mm}
{\rm{mod }}_{\sigma}(V_A, V_B)\leq C_1\, {\rm{Mod }}(A, B).
\end{equation}
\end{lemma}
\emph{Proof.} Let $\rho\co S(T)\rightarrow [0,\infty]$ be an
admissible Borel function for the pair $(A,B)$, i.e.,
$$
\int_{\gamma}\rho\, ds\geq 1,
$$
for every rectifiable curve $\gamma$ that connects $A$ and $B$. We
consider the mass distribution on $\sigma$ defined by
$$
m(v)=\bigg(\int_{St_2(v)}\rho^2d\mu\bigg)^{1/2}.
$$
To prove~(\ref{E:Mm}) we need to establish a mass bound and verify
admissibility. The mass bound is
$$
\begin{aligned}
\sum_{v\in V}m(v)^2 &\leq \sum_{v\in V}\bigg(\sum_{u\co
\delta(u,v)<2}\int_{U_u}\rho^2 d\mu \bigg)\\ &\lesssim \sum_{v\in
V}\int_{U_v}\rho^2 d\mu \lesssim \int_{S(T)}\rho^2 d\mu,
\end{aligned}
$$
where the constants understood depend only on $k$.

To show admissibility, we let $v_1, v_2,\dots, v_k$ be vertices of
a chain in $\sigma$ that connect $V_A$ and $V_B$. Then
$U_{v_1}\cap A\neq\emptyset,\ U_{v_k}\cap B\neq\emptyset$, and
$U_{v_{i-1}}\cap U_{v_i}\neq\emptyset$. We set $\lambda_1=A,\
\lambda_{k+1}=B$, and for $i=2,\dots, k$, let $\lambda_i$ be a
square in the square subdivision $\sigma$ with two of the vertices
being $v_{i-1}$ and $v_i$. Then for $i=2,\dots, k$, we have
$\lambda_i\in U_{v_{i-1}}\cap U_{v_i}$, and ${\rm{diam }}
\lambda_i=\sqrt{2}$. Also, since $A$ and $B$ are not contained in
any $St_2(v)$, there exists an absolute constant $c_0>0$ such that
${\rm{diam }} A\geq c_0$ and ${\rm{diam }} B\geq c_0$. Using
Lemma~\ref{L:Ad} we can inductively find rectifiable curves
$\eta_1,\dots, \eta_k$, satisfying the condition
$$
\int_{\eta_i}\rho\, ds\leq C_0 m(v_i),
$$
and such that $\eta_i$ connects
$\lambda_1\cup\eta_1\cup\dots\cup\eta_{i-1}$ and $\lambda_{i+1}$.
The constant $C_0$ depends only on $c_0$ and $k$. The union
$\eta_1\cup\dots\cup\eta_k$ contains a rectifiable curve $\eta$
connecting $A$ and $B$, and having the property
$$
1\leq\int_{\eta}\rho\, ds\leq C_0\sum_{i=1}^k m(v_i).
$$
Thus $C_0 m$ is an admissible mass distribution for the pair
$(V_A, V_B)$, and the proof is complete. \qed

\section{Proof of Theorem~\ref{T:main}}\label{S:Pr}
Let $M(r),\ r\geq0$, be an arbitrary $\R$-valued function, and
$L(\epsilon)$ be a function that satisfies the conditions of
Lemma~\ref{L:KeyL}, and such that $L(4\pi C_1/\log r)\geq M(r)$,
where $C_1$ is the constant from Lemma~\ref{L:Mm} when $k=3$. Let
$T$ be the subtree of $T_3$ given by Lemma~\ref{L:KeyL}. Then
$S(T)$ is parabolic, and let $g$ be a uniformizing map from $\C$
to $S(T)$. Let $A_{r'}$ and $B_r$ be the images under $g$ of
circles $\mathcal{C}_{r'}$ and $\mathcal{C}_r$ centered at 0 of
radii $r'$ and $r$, respectively, $1<r'<r$. We choose $r'$ such
that $A_{r'}$ is not contained in any set $St_2(v),\ v\in
V_{\sigma}$. Using Lemma~\ref{L:Mm} and the conformal invariance
of {\rm{Mod}}, we obtain that
$$
{\rm{mod }}_{\sigma}(V_{A_{r'}}, V_{B_r}) \leq C_1 {\rm{Mod }}
(\mathcal{C}_{r'}, \mathcal{C}_r)<\frac{4\pi C_1}{\log r},\ \
r\geq r_0=(r')^3.
$$
Since $T$ is a subgraph of $\sigma$, from monotonicity we have
$$
{\rm{mod }}_T(V_{A_{r'}}, V_{B_r}) <\frac{4\pi C_1}{\log r},\ \
r\geq r_0.
$$
If $D$ is the domain in $T$ which is the connected component of
$V_T\setminus V_{B_r}$ containing $v_0$, then ${\rm{mod
}}_T(\{v_0\},\dee D)< 4\pi C_1/\log r$. Therefore, by
Lemma~\ref{L:KeyL}, $|D|>L(4\pi C_1/\log r)\geq M(r)$. The proof
is complete.   \qed

\end{document}